\documentclass{amsart}

\usepackage{amscd,amssymb,epsfig}
\usepackage[all]{xy}

\def\comment#1{{\sf[#1]}}

\def\Z{{\mathbb Z}}
\def\Q{{\mathbb Q}}
\def\R{{\mathbb R}}
\def\C{{\mathbb C}}
\def\P{{\mathbb P}}
\def\F{{\mathbb F}}

\def\A{{\mathcal A}}
\def\H{{\mathcal H}}
\def\J{{\mathcal J}}
\def\L{{\mathcal L}}
\def\M{{\mathcal M}}
\def\T{{\mathcal T}}
\def\X{{\mathcal X}}

\def\G{{\Gamma}}
\def\S{{\Sigma}}
\def\h{{\mathfrak h}}

\def\Abar{\overline{\A}}
\def\Dtilde{\widetilde{D}}
\def\Fbar{\overline{\F}}
\def\Jbar{\overline{\J}}
\def\Mbar{\overline{\M}}
\def\Mtilde{\widetilde{\M}}

\def\That{{\widehat{T}}}
\def\Xbar{{\overline{X}}}
\def\Wtilde{\widetilde{W}}

\def\red{{\mathrm{red}}}
\def\sing{{\mathrm{sing}}}

\def\dot{{\bullet}}
\def\bs{\backslash}
\def\sing{{\mathrm{sing}}}

\newcommand\spec{\operatorname{Spec}}
\newcommand\Jac{\operatorname{Jac}}
\newcommand\Pic{\operatorname{Pic}}
\newcommand\Hom{\operatorname{Hom}}
\newcommand\Ext{\operatorname{Ext}}

\newcommand\Ind{\operatorname{Ind}}

\newtheorem{theorem}{Theorem}
\newtheorem{lemma}[theorem]{Lemma}
\newtheorem{proposition}[theorem]{Proposition}
\newtheorem{corollary}[theorem]{Corollary}

\theoremstyle{definition}

\theoremstyle{remark}

\begin{document}

\title[The Cohomology of the Space of Abelian $3$-folds]{The Rational
Cohomology Ring of the Moduli Space of Abelian $3$-folds}

\author{Richard Hain}
\address{Department of Mathematics\\ Duke University\\
Durham, NC 27708-0320}
\email{hain@math.duke.edu}

\thanks{Supported in part by grants from the National Science Foundation.}



\maketitle

\section{Introduction}

Suppose that $R$ is a commutative ring with unit.
Denote by $Sp_n(R)$ the group of automorphisms of $R^{2n}$ that preserve
the unimodular alternating form given by the matrix
$$
\begin{pmatrix}
0 & I_n \cr -I_n & 0
\end{pmatrix}
$$
In this note we compute the rational cohomology ring of $Sp_3(\Z)$, or
equivalently, of $\A_3$, the moduli space of principally polarized
abelian 3-folds.

Denote the $\Q$-Hodge structure of dimension 1 and weight $-2n$ by
$\Q(n)$. (For those not interested in Hodge theory, just interpret this
as one copy of $\Q$.) Denote by $\lambda$ the first Chern class in
$H^\dot(\A_3;\Q)$ of the Hodge bundle $\det \pi_\ast \Omega^1$ associated
to the projection $\pi$ of the universal abelian 3-fold to $\A_3$.

\begin{theorem}
\label{main}
The cohomology groups of $\A_3$ are given by
$$
H^j(\A_3;\Q) \cong H^j(Sp_3(\Z);\Q) \cong
\begin{cases}
\Q & j = 0; \cr
\Q(-1) & j = 2; \cr
\Q(-2) & j = 4; \cr
E & j = 6; \cr
0 & \text{otherwise,}
\end{cases}
$$
where $E$ is a two-dimensional mixed Hodge structure which is an extension
$$
0 \to \Q(-3) \to E \to \Q(-6) \to 0.
$$
The ring structure is determined by the condition that $\lambda^3 \neq 0$.
\end{theorem}

I do not know whether the mixed Hodge structure (MHS) $E$ on $H^6$ is split.
Since $\A_3$ is a smooth stack over $\spec \Z$, I expect it to be a
multiple (possibly trivial) of the class
$$
\zeta(3) \in \C/i\pi^3\Q \cong \Ext^1_\H(\Q,\Q(3))
$$
given by the value of the Riemann zeta function at 3. Determining this class
would be interesting.

As a corollary, we deduce the rational cohomology of $\Abar_3$, the Satake
compactification of $\A_3$.

\begin{theorem}
\label{satake_3}
The rational cohomology ring of $\Abar_3$ is given by
$$
H^j(\A_3;\Q) \cong
\begin{cases}
\Q(-n) & j = 2n, n \in \{ 0, 1, 2, 4, 5, 6\}; \cr
B & j = 6; \cr
0 & \text{otherwise,}
\end{cases}
$$
where $B$ is a 3-dimensional mixed Hodge structure which is an extension
$$
0 \to \Q(0) \to B \to \Q(-3)^2 \to 0.
$$
The ring structure is determined by the condition that the cohomology ring
contain the graded ring $\Q[\lambda]/(\lambda^7)$, where $\lambda$ has
degree 2 and type $(1,1)$.
\end{theorem}

Along the way, we compute the rational cohomology of $\Abar_2$, the Satake
compactification of $\A_2$ as well.

\begin{proposition}
\label{satake_2}
The rational cohomology ring of $\Abar_2$ is given by
$$
H^\dot(\Abar_2;\Q) \cong
\Q[\lambda]/(\lambda^4).
$$
where $\lambda$ is the first Chern class of the Hodge bundle
\end{proposition}

Denote the moduli space of smooth projective curves over the complex
numbers by $\M_g$. Another consequence of the proof is the surjectivity of
the homomorphism
$$
H^\dot(\A_3;\Q) \to H^\dot(\M_3;\Q)
$$
induced by the period mapping $\M_3 \to \A_3$.

The computation of the rational cohomology of $\A_1$ is classical\footnote{It
is trivial except in degree 0.} and follows from the fact that the quotient of
the upper half plane by $Sp_1(\Z)=SL_2(\Z)$ is a copy of the affine line. The
computation of the rational cohomology of $\A_2$ is (essentially) due to Igusa
\cite{igusa}.\footnote{It is 1-dimensional in degrees 0 and 2, and trivial
elsewhere.} Brownstein and Lee \cite{brown-lee,brown-lee:2} have computed the
integral cohomology of $Sp_2(\Z)$.

Suppose that $g\ge 2$. Recall that that the {\it mapping class group} $\G_g$ in
genus $g$ is the group of isotopy classes of orientation preserving
diffeomorphisms of a closed, oriented surface  $S$ of genus $g$. Its rational
cohomology is isomorphic to that of $\M_g$. The {\it Torelli group} $T_g$ is
defined to be the kernel of the natural homomorphism
\begin{equation}
\label{canon_homom}
\G_g \to Sp(H_1(S;\Z))
\end{equation}
where $Sp$ denotes the symplectic group, and where $H_1(S;\Z)$ is regarded
as a symplectic module via its intersection form. Choosing a symplectic
basis of $H_1(S;\Z)$ gives an isomorphism $Sp_g(\Z) \cong Sp(H_1(S;\Z))$.
One then obtains the well known extension
\begin{equation}
\label{xxx}
1 \to T_g \to \G_g \to Sp_g(\Z) \to 1.
\end{equation}
The {\it extended Torelli group} $\That_g$ is the preimage of the center
$\{\pm I\}$ of $Sp(H_1(S;\Z))$ under (\ref{canon_homom}).
Equivalently, it is the group of isotopy classes of diffeomorphisms of $S$
that act as $\pm I$ on $H_1(S)$. One has the extensions
$$
1 \to T_g \to \That_g \to \{\pm I\} \to 1
$$
and
\begin{equation}
\label{extension}
1 \to \That_g \to \G_g \to PSp_g(\Z) \to 1
\end{equation}
where $PSp_g(\Z)$ denotes the integral projective symplectic group
$Sp_g(\Z)/\{\pm I\}$.

Our approach to computing the cohomology of $Sp_3(\Z)$ is to analyse the
spectral sequence of the extension (\ref{extension}). This entails knowing
the cohomology of $\G_3$ (or equivalently, $\M_3$) and of $\That_3$.
Looijenga \cite{looijenga} computed the cohomology of $\M_3$ using
the theory of Del-Pezzo surfaces:
$$
H^j(\M_3;\Q) \cong
\begin{cases}
\Q & j = 0; \cr
\Q(-1) & j = 1; \cr
\Q(-6) & j = 6; \cr
0 & \text{otherwise.}
\end{cases}
$$

The cohomology of $\That_3$ is computed in Section~\ref{coho} using the
stratified Morse Theory of Goresky and MacPherson \cite{gor-macp}. We use their
theory of non-proper Morse functions (see their Part II, Chapter~10) applied to
the square of the distance to a point function restricted to the jacobian
locus. This is an elaboration of a trick of Geoff Mess \cite{mess} which he
used to show that the Torelli group in genus 2 is free of countable rank, and
which was also used by Johnson and Millson (cf.\ \cite{mess}) to show that
Torelli space in genus 3 does not have the homotopy type of a finite complex.
Our use of Morse theory is analogous to Goresky and MacPherson's treatment
\cite[Part III]{gor-macp} of complements of affine subspaces of euclidean
spaces. \comment{Look at Mess's paper.}

It seems to be a curious fact that in low genus ($g=2,3$ so far), it
is easier to compute the rational cohomology groups of $\M_g$ than
of $\A_g$. This is perhaps a reflection of the richness of curve theory ---
that it is a more powerful tool for understanding the geometry of $\M_g$ than
the theory of abelian varieties is as a tool for understanding the geometry
of $\A_g$. It will be interesting to know if this trend persists when
$g \ge 4$, when $\M_g$ is no longer dense in $\A_g$. The cohomology of $\M_4$
is not yet know, nor does it seem tractable to compute the cohomology
of the extended Torelli group in genus 4.

\section{Preliminaries}

General references for this section are  \cite[Chapt.~3]{griffiths-harris}
and \cite{hain-looijenga}.

A {\it polarized abelian variety} is a compact complex torus $A$ together
with a cohomology class
$$
\theta \in H^{1,1}(A) \cap H^2(A;\Z)
$$
whose translation invariant representative is positive. The corresponding
complex line bundle is ample. The polarization $\theta$ can be regarded
as a skew symmetric bilinear form on $H_1(A;\Z)$. The polarization is
{\it principal} if this form is unimodular.

Jacobians of curves are polarized by the intersection pairing on $H_1(C;\Z)
\cong H_1(\Jac C;\Z)$. This form is unimodular, and so jacobians are
canonically principally polarized abelian varieties.

A {\it framed} principally polarized abelian variety is a principally
polarized abelian variety $A$ together with a symplectic basis of
$H_1(A;\Z)$ with respect to the polarization $\theta$. 

Suppose that $g\ge 1$.
The maximal compact subgroup of $Sp_g(\R)$ is $U(g)$. The symmetric
space $Sp_g(\R)/U(g)$ is isomorphic to the rank $g$ Siegel upper half space
$$
\h_g = \left\{
\parbox{2.25in}{symmetric $g\times g$ complex matrices with
positive definite imaginary part}
\right\}
$$
(See, \cite{helgason}, for example.) It has dimension $g(g+1)/2$.
Taking a framed principally polarized abelian variety
$(A;a_1,\dots,a_g,b_1,\dots,b_g)$ to the corresponding period matrix
gives a bijection
\begin{equation}
\label{bijection}
\h_g \cong \left\{
\parbox{2.25in}{isomorphism classes of framed principally polarized abelian
varieties}
\right\}
\end{equation}
We will regard $\h_g$ as the (fine) moduli space of framed 
principally polarized abelian varieties of dimension $g$.

Changing symplectic bases gives a natural left action of $Sp_g(\Z)$ on the
moduli space of framed principally polarized abelian varieties. There is
clearly a natural left $Sp_g(\Z)$ action on $\h_g = Sp_g(\Z)/U(g)$. The
bijection (\ref{bijection}) is equivariant with respect to these actions.

The moduli space $\A_g$ of principally polarized abelian varieties of
dimension $g$ is the quotient $Sp_g(\Z)\bs \h_g$. Since $\h_g$ is contractible
and since $Sp_g(\Z)$ acts discontinuously and virtually freely on $\h_g$,
it follows that there is a natural isomorphism
$$
H^\dot(Sp_g(\Z);\Q) \cong H^\dot(\A_g;\Q).
$$
Taking a curve $C$ to its jacobian $\Jac C$ defines a morphism $\M_g\to \A_g$
which is called the {\it period mapping}.

Now suppose that $g\ge 2$. Denote Teichm\"uller space in genus $g$ by $\X_g$.
The mapping class group $\G_g$ acts properly discontinuously and virtually
freely on $\X_g$ with quotient $\M_g$. It follows that there is a natural
isomorphism
$$
H^\dot(\G_g;\Q) \cong H^\dot(\M_g;\Q).
$$

We shall need several moduli spaces that sit between $\X_g$ and $\M_g$. Denote
the quotient of $\X_g$ by $T_g$ by $\T_g$. This space is known as {\it Torelli
space}. Since $\X_g$ is contractible and $T_g$ is torsion free, $T_g$ acts
freely on $\X_g$ and Torelli space is an Eilenberg-MacLane space with
fundamental group $T_g$. Consequently,
$$
H^\dot(T_g;\Z) \cong H^\dot(\T_g;\Z).
$$

By a {\it framed Riemann surface} of genus $g$ we shall mean a compact
Riemann surface $C$ together with a symplectic basis $a_1,\dots,a_g, b_1,
\dots, b_g$ of $H_1(C;\Z)$ with respect to the intersection form. Torelli
space $\T_g$ is the moduli space of framed Riemann surfaces of genus $g$;
its points  correspond to isomorphism classes of framed, genus $g$ Riemann
surfaces. The symplectic group $Sp_g(\Z)$ acts on the framings in the natural
way; the quotient $Sp_g(\Z)\bs \T_g$ is $\M_g$.

Denote the locus in $\h_g$ consisting of jacobians of smooth curves by
$\J_g$. (Note that this is not closed in $\A_g$.) The period mapping
$\T_g \to \J_g$ is surjective by definition.
Since minus the identity is an automorphism of every polarized abelian variety
$A$,
$$
(A;a_1,\dots,b_g) \cong (A;-a_1,\dots, -b_g).
$$
But if $C$ is a genus $g$ curve, then
$$
(C;a_1,\dots,b_g) \cong (C;-a_1,\dots, -b_g)
$$
if and only if $C$ is hyperelliptic. It follows that, when $g\ge 3$, the
period mapping
$$
\T_g \to \J_g
$$
that takes a framed curve to its jacobian with the same framing is
surjective and 2:1 except along the hyperelliptic locus, where it is 1:1.
It follows that, when $g\ge 3$,
$$
\J_g = \That_g \bs \X_g \text{ and }\M_g = PSp_g(\Z) \bs \J_g.
$$

The following diagram shows the coverings and their Galois group
when $g\ge 3$:
\begin{equation}
\label{covers}
\xymatrix{
\X_g \ar[rr]_{T_g}\ar@/^1pc/[rrrr]^{\That_g}\ar@/^3pc/[rrrrrr]^{\G_g} &&
\T_g \ar[rr]_{\Z/2\Z}\ar[rrdd]\ar@/^1pc/[rrrr]^{Sp_g(\Z)} && \J_g
\ar[rr]_{PSp_g(\Z)}\ar@{_{(}->}[dd] && \M_g \ar[dd] \cr \cr
& & & &\h_g \ar[rr]^{PSp_g(\Z)} & & \A_g
}
\end{equation}

\begin{lemma}
\label{quot}
If $1 \to \Z/2\Z \to E \to G \to 1$
is a group extension, then the projection $E\to G$ induces an isomorphism
on homology and cohomology with 2-divisible coefficients. In particular
$$
H^\dot(PSp_g(\Z);\Z[1/2]) \to H^\dot(Sp_g(\Z);\Z[1/2])
$$
is an isomorphism.
\end{lemma}

\begin{proof}
This follows from the fact that
$$
H^j(\Z/2;\Z[1/2]) = 0 \qquad j > 0
$$
using the Hochschild-Serre spectral sequence of the group extension.
\end{proof}

\begin{proposition}
For all $g \ge 2$, there is a natural isomorphism
$$
H_\dot(\J_g;\Z[1/2]) \cong H_\dot(\That_g;\Z[1/2])
\cong H_\dot(T_g;\Z[1/2])^{\Z/2\Z}.
$$
There are similar isomorphisms for cohomology.
\end{proposition}

\begin{proof}

Since $T_g$ is torsion free, $T_g$ acts fixed point freely on Teichm\"uller
space, and $\T_g$ is a model of the classifying space of $T_g$. Recall that
if $X$ is a simplicial complex on which $\Z/2$ acts simplicially (but not
necessarily fixed point freely), then the map 
$$
p_\ast : H_\dot(X/(\Z/2);\Z[1/2]) \to  H_\dot(X;\Z[1/2])^{\Z/2}
$$
induced by the projection $p$ is an isomorphism, whose inverse is half the
pullback map $p^\ast$. Applying this twice gives isomorphisms
$$
H^\dot(\J_g;\Z[1/2]) \cong H^\dot(\T_g;\Z[1/2])^{\Z/2}
\cong H^\dot(T_g;\Z[1/2])^{\Z/2} \cong H^\dot(\That_g;\Z[1/2]).
$$
\end{proof}

A theta divisor of a principally polarized abelian variety $A$ is a divisor
$\Theta$ whose Poincar\'e dual is the polarization and which satisfies
$i^\ast \Theta = \Theta$, where $i:x \mapsto -x$. Any two such divisors
differ by translation by a point of order 2, and can be given as the
zero locus of a theta function associated to a period matrix of $A$.
A principally polarized abelian variety $A$ is {\it reducible} if it is
isomorphic (as a polarized variety) to the product of two proper abelian
subvarieties. If $A = A_1 \times A_2$, then any theta divisor of $A$
is reducible:
$$
\Theta_A =
\big(\Theta_{A_1} \times A_2\big) \cup \big(A_1 \times \Theta_{A_2}\big).
$$
Denote the locus of reducible abelian varieties in $\A_g$ by $\A_g^\red$
and in $\h_g$ by $\h_g^\red$. Elements of $\h_g^\red$ are precisely those
period matrices $\Omega$ that can be written as a direct sum of two 
smaller period matrices.

\begin{proposition}
Denote the closure of $\J_g$ in $\h_g$ by $\Jbar_g$.
If $g\ge 2$, then $\J_g = \Jbar_g - (\Jbar_g \cap \h^\red_g)$.
\end{proposition}

\begin{proof}
The period mapping $\M_g \to \A_g$ extends to a morphism $\Mbar_g \to \Abar_g$
from the Deligne-Mumford compactification of $\M_g$ to the Satake
compactification of $\A_g$. The inverse image of the boundary $\Abar_g - \A_g$
of $\Abar_g$ is the boundary divisor $\Delta_0$ of $\Mbar_g$, whose generic
point is an irreducible stable curve of genus $g$ with one node. Denote the
moduli space of curves of compact type $\Mbar_g - \Delta_0$ by $\Mtilde_g$.
Since $\Mbar_g$ is complete, it follows that the period mapping $\Mtilde_g \to
\A_g$ is proper and therefore has closed image in $\A_g$. Since $\M_g$ is dense
in $\Mbar_g$ and has image $Sp_g(\Z)\bs\J_g$ under the period mapping, it
follows that the image of $\Mtilde_g$ in $\A_g$ is $Sp_g(\Z)\bs\Jbar_g$.

Recall that the theta divisor $\Theta_C \subset \Jac C$ of a smooth genus $g$
curve $C$ is (up to a translate by a point of order 2) the image of mapping
$$
C^{g-1} \to \Pic^{g-1} C \to \Jac C
$$
that takes $(x_1,\dots,x_{g-1})$ to $x_1 + \dots + x_{g-1} - \alpha$, where
$\alpha$ is a square root of the canonical bundle of $C$.
(See, for example, \cite[p.~338]{griffiths-harris}.) It follows that
$\Theta_C$ is irreducible. On the other hand, if $C$ is a reducible,
stable, curve of compact type, its jacobian is the product of the
components of its irreducible components, and is therefore reducible. The
result follows.
\end{proof}

\begin{corollary}
\label{j3}
We have $\J_3 = \h_3 - \h_3^\red$.
\end{corollary}

\begin{proof}
Both $\M_3$ and $\A_3$ have dimension 6. Since the period mapping
$\M_3 \to \A_3$ is generically of maximal rank, $\Mtilde_3 \to \A_3$
is surjective. This implies that $\Jbar_3 = \h_3$, from which the result
follows.
\end{proof}

\section{The Homology of $\J_3$ and $\That_3$}
\label{coho}

Denote the singular locus of an analytic variety $Z$ by $Z^\sing$. We shall
compute the homology of $\J_3$ by applying stratified Morse theory to the
the stratification
\begin{equation}
\label{strat}
\h_3 \supseteq \h_3^\red \supset \h_3^{\red,\sing}
\end{equation}
of $\h_3$. The top stratum is, by Corollary~\ref{j3}, $\J_3$. Note that there
are natural inclusions
$$
\h_1 \times \h_2 \hookrightarrow \h_3^\red \text{ and }
\h_1 \times \h_1 \times \h_1 \hookrightarrow \h_3^{\red,\sing}
$$
defined by
$$
(\tau,\Omega) \mapsto
\begin{pmatrix}
\tau & 0 \cr 0 & \Omega
\end{pmatrix}
\text{ and }
(\tau_1,\tau_2,\tau_3) \mapsto
\begin{pmatrix}
\tau_1 & 0 & 0 \cr 0 & \tau_2 & 0 \cr 0 & 0 & \tau_3
\end{pmatrix}
$$
respectively. Note that the image of $\h_1 \times \h_2$ is stabilized in
$Sp_3(\R)$ by $SL_2(\R)\times Sp_2(\R)$, and the image of $(\h_1)^3$ by
$\Sigma_3 \ltimes SL_2(\R)^3$,
where $\Sigma_3$ is identified with the subgroup
$$
\{a_j \mapsto a_{\sigma(j)} \text{ and } b_j \mapsto b_{\sigma(j)} :
\sigma \text{ is a permutation of }\{1,2,3\}\}
$$
of $Sp_3(\R)$. Here $a_1,\dots, b_6$ is the distinguished framing of the
first homology of the corresponding abelian variety.

\begin{proposition}
\label{strata}
The stratification (\ref{strat}) satisfies Whitney's conditions (A) and (B)
(cf.\ \cite[p.~37]{gor-macp}). Moreover
\begin{equation}
\label{red}
\h_3^\red = \bigcup_{g\in Sp_3(\Z)} g(\h_1 \times \h_2)
= \bigcup_{g\in Sp_3(\Z)/\left(SL_2(\Z)\times Sp_2(\Z)\right)}
g(\h_1 \times \h_2)
\end{equation}
and
\begin{equation}
\label{sing}
\h_3^{\red,\sing} = \bigcup_{g\in Sp_3(\Z)} g(\h_1 \times \h_1 \times \h_1)
= \coprod_{g\in Sp_3(\Z)/\left(\Sigma_3\ltimes (SL_2(\Z))^3\right)}
g(\h_1 \times \h_1 \times \h_1)
\end{equation}
In particular, $\h_3^\red$ is a locally finite union of totally geodesic
complex submanifolds of $\h_3$ of complex codimension $2$ and
$\h_3^{\sing,\red}$ is a countable disjoint union of totally geodesic complex
submanifolds of $\h_3$ of  complex codimension $3$.
\end{proposition}

\begin{proof}
Since every reducible abelian variety is the product (as polarized varieties)
of an elliptic curve and an abelian surface, $\A_1 \times \A_2 \to \A_3^\red$
is surjective. Lifting to $\h_3$, this implies that the $Sp_3(\Z)$ acts
transitively on the components of $\h_3^\red$, and that $\h_3^\red$ is the
$Sp_3(\Z)$-orbit of $\h_1 \times \h_2$ in $\h_3$. Since $SL_2(\R) \times
Sp_2(\R)$ acts transitively on $\h_1 \times \h_2$, the stabilizer in $Sp_3(\Z)$
of $\h_1 \times \h_2$ is $SL_2(\Z) \times Sp_2(\Z)$. Assertion (\ref{sing})
follows.

The components of $\h_3^\red$ are smooth, so the $\h_3^{\red,\sing}$ is the
locus where two or more components of $\h_3^\red$ intersect. This is precisely
the preimage of the locus in $\A_3$ of products of 3 elliptic curves. Since
this locus is irreducible (it is the image of $(\A_1)^3 \to \A_3$),
$\h_3^{\red,\sing}$ is the $Sp_3(\Z)$-orbit of $(\h_1)^3$. By the semi-simplicity
of polarized abelian varieties, there is a unique way to decompose an element
of $\h_3^{\red,\sing}$ as a product of three elliptic curves. This, and the fact
that $SL_2(\R)^3$ acts transitively on $(\h_1)^3$, imply the stabilizer in
$Sp_3(\Z)$ of $(\h_1)^3$ is $\S_3 \ltimes SL_2(\Z)^3$ and that each component of
$\h_3^{\red,\sing}$ is smooth. This proves (\ref{red}).

Whitney's condition (A) is automatic as each component of $\h_3^{\red,\sing}$
is a homogeneous submanifold of the closure of each component of $\h_3^\red$.
Condition (B) is well known to be a consequence of condition (A).
\end{proof}

Integrating the Riemannian metric along geodesics gives an $Sp_g(\R)$-invariant
distance function $d$ on $X$. For a point $p \in X$, let $D_p : X \to \R$ be
the square of the distance to $p$:
$$
D_p(x) = d(x,p)^2.
$$

The following result can be proved, either by appealing to
\cite[I.2.2.3]{gor-macp} or by an elementary and direct argument.

\begin{proposition}
There is an open dense subset $U$ of $X$ such that for all $p \in U$,
$D_p : X \to \R$ is a Morse function in the sense of Goresky and MacPherson
\cite[p.~52]{gor-macp}, all of whose critical points are ``nondepraved.'' \qed
\end{proposition}

Since each stratum is a union of totally geodesic subspaces, and since
the symmetric space metric is complete with non-positive curvature it follows
that there is a unique critical point on each component of each stratum.
Since the Morse data for each critical point is a product of the normal
and tangential Morse data \cite[p.~61]{gor-macp}, we only
need compute the normal Morse data at each critical point. There are two
types of these: those that lie on a translate of $\h_1\times \h_2$ and
those that lie on a translate of $(\h_1)^3$. The Morse data at each 
depends only on its type.

\begin{proposition}
If $x \in \h_3^\red - \h_3^{\red,\sing}$ is a critical point of $D_p$,
then the normal Morse data at $x$ is homotopy equivalent to $(S^3,*)$.
In particular, $D_p$ is perfect at such critical points.
\end{proposition}

\begin{proof}
Since the normal slice at a smooth point of $\h_3^\red$ is a complex 2-ball,
the normal Morse data at $x$ is $(S^3,*)$. Since $H_\dot(S^3) \to H_\dot(S^3,*)$
is surjective, $D_p$ is perfect at $x$.
\end{proof}

\begin{lemma}
At each point of $\h_3^{\red,\sing}$, there is a normal slice with coordinates
$(z_1,z_2,z_3)$ such that $\h_3^\red$ has three components with equations
$$
z_2 = z_3 = 0,\quad z_1 = z_3 = 0,\quad z_1 = z_2 = 0.
$$
\end{lemma}

\begin{proof}
There are three obvious ways to deform the product $A=E \times E' \times E''$
of three elliptic curves, preserving the polarization, into $\h_3^\red$. Namely,
one can deform one of the elliptic curves in $\h_1$, and deform the product of
the other two into $\h_2 - \h_2^\red$. The semi-simplicity of abelian varieties
implies that each component of $\h_3^{\red,\sing}$ is smooth and there are no
other ways to deform $A$ into $\h_3^\red$. It follows that 3 components of
$\h_3^\red$ intersect at each point of $\h_3^{\red,\sing}$.

\begin{figure}[!ht]
\epsfig{file=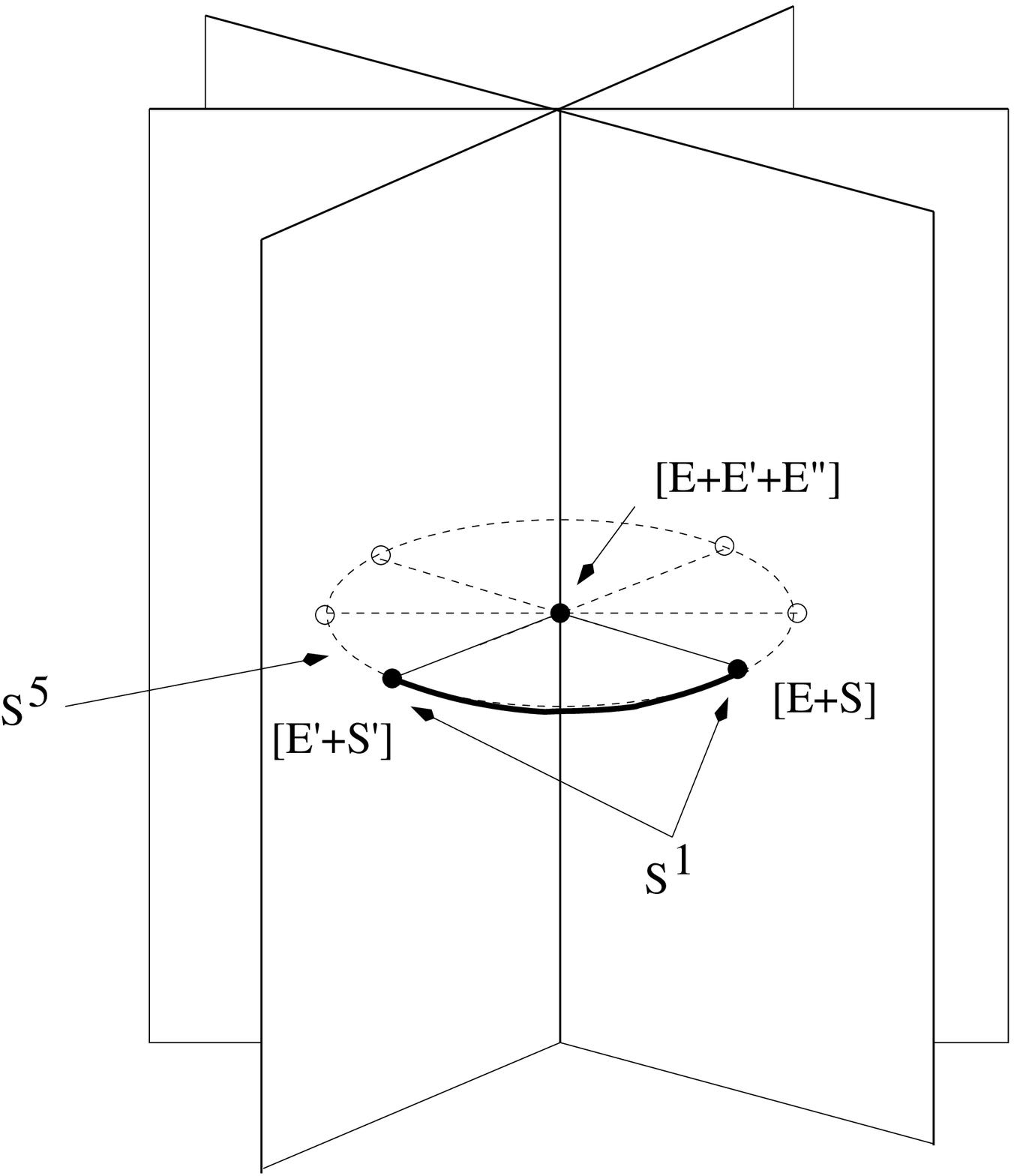, width=3.0in}
\caption{Local structure of $\h_3^\red$ near $\h_3^{\red,\sing}$}
\end{figure}

Since each component of $\h_3^{\red,\sing}$ is locally homogeneous, and since
all are conjugate under the action of $Sp_3(\Z)$, to determine the local
structure of $\h_3^\red$ near $\h_3^{\red,\sing}$, it suffices to write down
local equations along the component $(\h_1)^3$ of $\h_3^\red$. Here we can
use the coordinates
$$
\begin{pmatrix}
\tau_1 & z_3 & z_2 \cr z_3 & \tau_2 & z_1 \cr z_2 & z_1 & \tau_3
\end{pmatrix}
$$
where each $\tau_j \in \h_1$, and $(z_1,z_2,z_3)$ lies in a neighbourhood of
the origin in $\C^3$ small enough to guarantee that this matrix has positive
definite imaginary part.
In these coordinates, the three components of $\h_3^\red$ have equations:
$$
z_2 = z_3 = 0,\quad z_1 = z_3 = 0,\quad z_1 = z_2 = 0.
$$
\end{proof}

View $S^5$ as the unit sphere in $\C^3$. The intersection of each
coordinate axis with $S^5$ is a linearly imbedded $S^1$. Together these
give an imbedding
$$
k : S^1 \amalg S^1 \amalg S^1 \hookrightarrow S^5
$$
where $\amalg$ denotes disjoint union.
The boundary of a small tubular neighbourhood of each $S^1$ is a trivial
$S^3$ bundle over $S^1$. Taking one fiber of each and connecting them to
an arbitrary point of $S^5 - \big(S^1 \cup S^1 \cup S^1\big)$ gives an
imbedding 
$$
i : S^3 \vee S^3 \vee S^3 \hookrightarrow S^5.
$$

\begin{proposition}
If $x \in \h_3^{\red,\sing}$ is a critical point of $D_p$,
then the normal Morse data at $x$ is homotopy equivalent to
$$
\big(S^5 - k\big(S^1\amalg S^1 \amalg S^1\big),
i(\big( S^3 \vee S^3 \vee S^3\big) \big).
$$
\end{proposition}

\begin{proof}
By the previous Lemma, the normal Morse data is the same as that for the
distance squared to a point function for
$$
\C^3 - \text{ the union of the 3 coordinate axes. }
$$
This is computed using \cite[I.3.11.2)]{gor-macp} as explained in
\cite[III.3.3)]{gor-macp}.
\end{proof}

To compute the homology of $\J_3$, we need to compute the homology of each kind
of Morse data. Denote the quotient of $\Z^3$ by the diagonal subgroup by $V$.
It is  isomorphic to $\Z^2$.

\begin{proposition}
The homology of the Morse data at a critical point of $\h_3^{\red,\sing}$ is
$$
H_j\big(S^5 - k\big(S^1\amalg S^1 \amalg S^1\big),
i\big( S^3 \vee S^3 \vee S^3\big);\Z \big) =
\begin{cases}
V & j = 4; \cr
0 & \text{otherwise}.
\end{cases}
$$
The homology is generated by the boundaries of tubular neighbourhoods of
the three imbedded $S^1$s, which are subject to the relation that their sum
is zero. In particular, $D_p$ is perfect at such critical points.
\end{proposition}

\begin{proof}
The computation is elementary. The Morse function is perfect at such critical
points because the relative homology is generated by absolute cycles.
\end{proof}

The natural action of the symmetric group $\S_3$ on $\Z^3$ (by permuting
the coordinates) preserves the diagonal subgroup and therefore descends
to an action on $V$. We view $V$ as a $\S_3\ltimes SL_2(\Z)^3$-module
via the projection $\S_3\ltimes SL_2(\Z)^3 \to \S_3$.

Recall that if $R$ is a commutative ring, $K$ a subgroup of $G$, and $M$
a $RK$ module, then the $G$-module induced from $M$ is defined by
$$
\Ind_K^G M := RG\otimes_{RK} M.
$$

\begin{theorem}
We have,
$$
H_j(\J_3;\Z) \cong
\begin{cases}
\Z & j = 0; \cr
\Ind_{SL_2(\Z)\times Sp_2(\Z)}^{Sp_3(\Z)} \Z & j = 3;\cr
\Ind_{\S_3 \ltimes SL_2(\Z)^3}^{Sp_3(\Z)} V & j = 4; \cr
0 & \text{otherwise.}
\end{cases}
$$
\end{theorem}

\begin{proof}
Each critical point has trivial tangential Morse data, as the critical
points are minima on each stratum. It follows from \cite[I.3.7]{gor-macp}
that the Morse data at each critical point is homotopy equivalent to its
normal Morse data. Since the Morse function is perfect, the homology
of $\J_3$ is the sum of the relative homologies of the normal Morse data at
each critical point. The action of $Sp_3(\Z)$ follows from the description
of the $Sp_3(\Z)$ action on the strata in Proposition~\ref{strata}.
\end{proof}

For a subgroup $G$ of $Sp_g(\Z)$ that contains $-I$, define $PG$ to be
the subgroup $G/(\pm I)$ of $PSp_g(\Z)$.

\begin{corollary}
For all $j \ge 0$,
\begin{align*}
H_j(\That_3;\Z[1/2]) \cong H_k(T_3;\Z[1/2])^{\Z/2\Z} &\cong
\begin{cases}
\Z[1/2] & j = 0; \cr
\Ind_{SL_2(\Z)\times Sp_2(\Z)}^{Sp_3(\Z)} \Z[1/2] & j = 3;\cr
\Ind_{\S_3 \ltimes SL_2(\Z)^3}^{Sp_3(\Z)} V\otimes\Z[1/2] & j = 4; \cr
0 & \text{otherwise;}
\end{cases}
\cr
&\cong
\begin{cases}
\Z[1/2] & j = 0; \cr
\Ind_{P(SL_2(\Z)\times Sp_2(\Z))}^{PSp_3(\Z)} \Z[1/2] & j = 3;\cr
\Ind_{P(\S_3 \ltimes SL_2(\Z)^3)}^{PSp_3(\Z)} V\otimes\Z[1/2] & j = 4; \cr
0 & \text{otherwise.}
\end{cases}
\end{align*}
\end{corollary}

\section{The Spectral Sequence}

We shall compute the homology spectral sequence
\begin{equation}
\label{ss}
H_i(PSp_3(\Z);H_j(\That_3;\Q)) \implies H_{i+j}(\G_3;\Q).
\end{equation}

Thanks to Shapiro's Lemma (see, for example, \cite{brown})
the $E^2$ term of the spectral sequence can be computed.

\begin{lemma}
We have
$$
H_i(PSp_3(\Z);H_j(\That_3;\Q)) \cong
\begin{cases}
\Q & i = 0, 2 \text{ and } j=3; \cr
0 &  j > 0 \text{ and } j\neq 3.
\end{cases}
$$
\end{lemma}

\begin{proof}

Applying Shapiro's Lemma, Lemma~\ref{quot}, the Kunneth Theorem, and the fact
that the rational homology of $SL_2(\Z)$ is that of a point, we have:
\begin{align*}
H_i\big(PSp_3(\Z);\Ind_{P(SL_2(\Z)\times Sp_2(\Z))}^{PSp_3(\Z)}\Q\big)
&\cong H_i(P(SL_2(\Z)\times Sp_2(\Z));\Q) \cr
&\cong H_i(SL_2(\Z)\times Sp_2(\Z);\Q)\cr
&\cong H_i(Sp_2(\Z);\Q).
\end{align*}
This is $\Q$ in when $i=0,2$ and 0 otherwise by Igusa's computation.

Let $V_\Q = V\otimes_\Z \Q$. This is the unique 2-dimensional irreducible
represention of $\S_3$. Since $V_\Q$ is divisible and has no coinvariants,
$H_\dot(\S_3;V_\Q)$ vanishes in all degrees. Arguing as above, we have
\begin{align*}
H_i\big(PSp_3(\Z);\Ind_{P(\S_3 \ltimes SL_2(\Z)^3)}^{PSp_3(\Z)}V_\Q)
&\cong H_i(P(\S_3 \ltimes SL_2(\Z)^3);V_\Q) \cr
&\cong H_i(\S_3 \ltimes SL_2(\Z)^3;V_\Q) \cr
&\cong H_i(\S_3;V_\Q) \cr
&\cong 0.
\end{align*}
\end{proof}

\begin{proposition}
\label{spec_seq}
For $2 \le r \le 4$, the $E^r$-term of the spectral sequence (\ref{ss}) is
\begin{center}
\begin{tabular}{c|cccccccc}	
$4$ & $0$ & $0$ & $0$ & $0$ & $0$ & $0$ & $0$ & $0$ \cr
$3$ & $\Q$ & $0$ & $\Q$ & $0$ & $0$ & $0$ & $0$ & $0$ \cr
$2$ & $0$ & $0$ & $0$ & $0$ & $0$ & $0$ & $0$ & $0$ \cr
$1$ & $0$ & $0$ & $0$ & $0$ & $0$ & $0$ & $0$ & $0$ \cr
$0$ & $\Q$ & $0$ & $\Q$ & $0$ & $\Q$ & $0$ & $\Q^2$ & $0$ \cr
\hline 
$\mathrm{deg}$  & $0$ & $1$ & $2$ & $3$ & $4$ & $5$ & $6$ & $7$ \cr
\end{tabular}
\end{center}
(All terms not shown are zero.) In addition, the differentials
$$
d^4 : E^4_{4,0} \to E^4_{0,3} \text{ and } d^4 : E^4_{6,0} \to E^4_{2,3}
$$
are both surjective, and $E^5 = E^\infty$.
\end{proposition}

\begin{proof}
The computation of $E^2_{s,t}$ for $t>0$ follows from the previous lemma.
It implies that $d^r = 0$ when $2\le r < 4$ and $E^5 = E^\infty$.
By Looijenga's computation of the rational homology of $\G_3$ we know that
$H_3(\G_3;\Q)$ and $H_5(\G_3;\Q)$ both vanish. This implies that the
differentials
$$
d^4 : E^4_{4,0} \to E^4_{0,3} \text{ and } d^4 : E^4_{6,0} \to E^4_{2,3}
$$
must be surjective. Since $H_4(\G_3;\Q) = 0$, and since $H_6(\G_3;\Q)$
is one dimensional, the first of these is an isomorphism and the second has
one dimensional kernel. The result follows. 
\end{proof}

\begin{proof}[Proof of Theorem~\ref{main}]
The computation of the rational homology (and therefore the rational
cohomology) of $\A_3$ follows from Proposition~\ref{spec_seq}.

Denote the first Chern class in $H^2(\A_3;\Q)$ of the Hodge bundle by
$\lambda$. It is the class of an ample line bundle. A standard argument, that
uses the fact that the Satake compactification of $A_3$ has boundary of
codimension 3, shows that there is a complete surface in $A_3$. This implies
that $\lambda^2 \neq 0$ in $H^4(\A_g;\Q)$.

The proof that $\lambda^3$ does not vanish in $H^6(\A_3;\Q)$ is more subtle,
and is due to van der Geer \cite{vandergeer}. We will give a topological proof
of this fact in the next section. The key point in van der Geer's argument is
that there is a complete subvariety of the characteristic $p$ version
$\A_{3/\Fbar_p}$ of $\A_3$. This implies that $\lambda^3$ is not zero in
$H_{\text{\'et}}^6(\A_{3/\Fbar_p};\Q_\ell)$, where $p$ is a prime where $\A_3$
has good reduction and $\ell$ is a prime distinct from $p$. Standard comparison
theorems imply that this last group is isomorphic to $H^6(\A_3;\Q_\ell)$, which
gives the desired non-vanishing of $\lambda^3$.

The statement about weights follows as $\lambda$ is of type $(1,1)$. Since
$\lambda^3 \neq 0$, this implies that $H^4(\A_3;\Q)$ is generated by
$\lambda^2$ and has type $(2,2)$, and that $H^6(\A_3;\Q)$ contains a copy of
$\Q(-3)$ spanned by $\lambda^3$. On the other hand, the spectral sequence in
Proposition~\ref{spec_seq} implies that the restriction mapping
$$
H^6(\A_3;\Q) \to H^6(\M_3;\Q) \cong \Q(-6)
$$
is surjective, which completes the proof.
\end{proof}

\section{Cycles}

In this section, we give a topological proof that $\lambda^3$ is non-zero in
$H^6(\A_3;\Q)$. The approach is to construct a topological 6-cycle in $\A_3$
and then show that the value of $\lambda^3$ on it is non-zero.

Fix an imbedding of the Satake compactification $\Abar_3$ of $\A_3$ in some
projective space. Take a generic codimension 3 linear section of $\Abar_3$ that
avoids the boundary $\Abar_3^\red - \A_3^\red$, is transverse to $\A_3^\red$,
and intersects $\A_3^{\red,\sing}$ transversally. This section is a complete
curve $X$ in $\A_3^\red$, smooth (in the orbifold sense) away from its
intersection with $\A_3^{\red,\sing}$. At each point $x$ where it intersects
$\A_3^{\red,\sing}$, it has three branches. Set $X' =
X - (X\cap \A_3^{\red,\sing})$.

\begin{figure}[!ht]
\epsfig{file=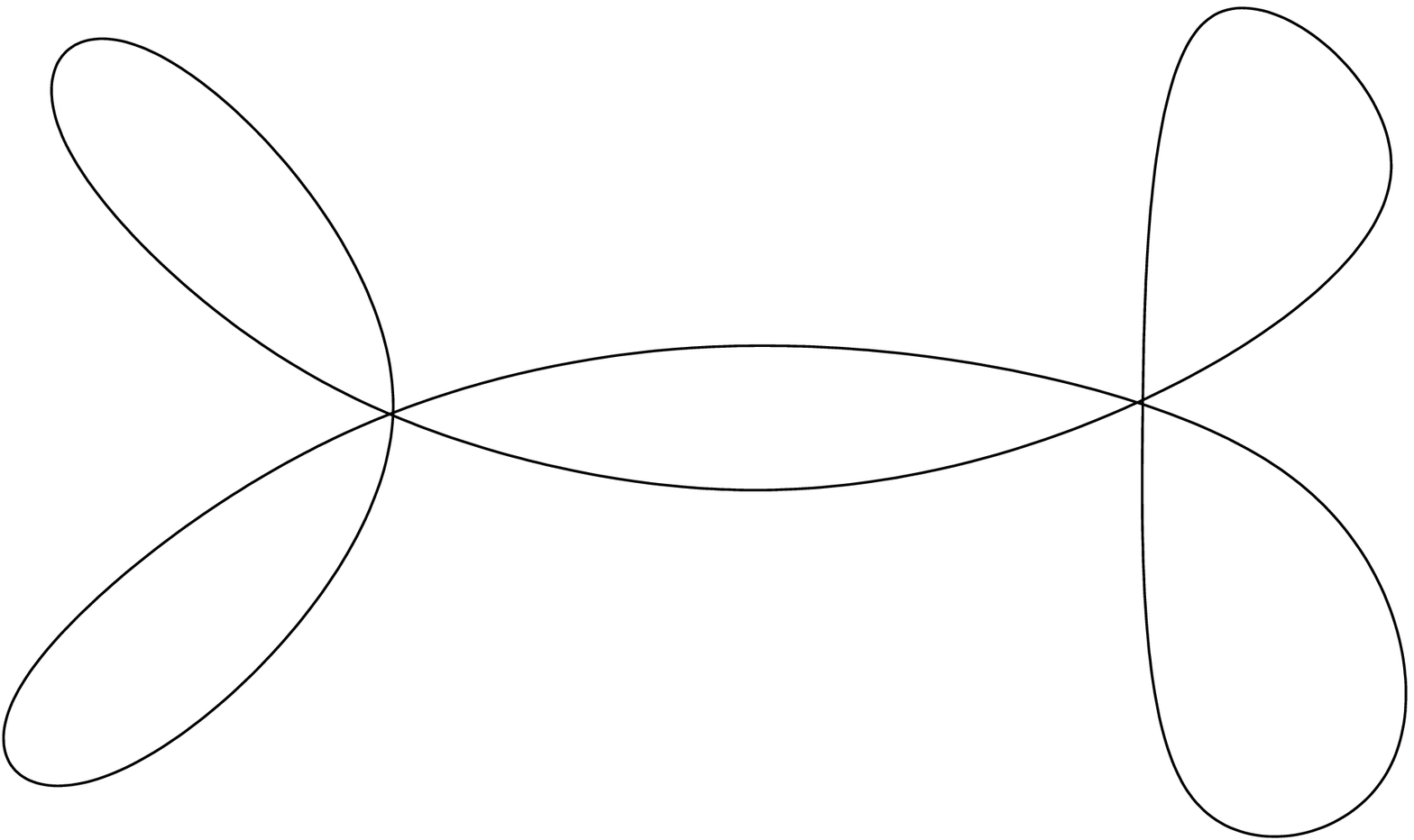, width=2.0in}
\caption{}
\end{figure}

The next step is to construct a 5-cycle in $\M_3$ which is an $S^3$ bundle over
$X$ away from the triple points, and where the three branches of this bundle
at each triple point $x$ are plumbed together using the normal Morse data at
$x$.

Denote the moduli space of principally polarized abelian 3-folds with a level
$\ell$ structure by $\A_3[\ell]$. This is the quotient of $\h_3$ by the level
$\ell$ subgroup $Sp_3(\Z)[\ell]$ of $Sp_3(\Z)$, and is smooth when $\ell \ge
3$. Fix an $\ell \ge 3$. The symmetric space metric on $\h_3$ descends to
$\A_3[\ell]$. With the help of the metric, the normal bundle of any stratum can
be viewed as a subbundle of the tangent bundle of $\A_3[\ell]$. Using the
exponential mapping, we can identify a neighbourhood of the zero section of the
normal bundle as being imbedded in $\A_3[\ell]$. Let $Y$ and $Y'$ be the
inverse images of $X$ and $X'$ in $\A_3[\ell]$.

Choose a positive real number $\epsilon$ such that the exponential mapping is an
imbedding on the $\epsilon/4$-ball $B$ of the normal bundle of
$\A_3^\red[\ell]$ restricted to $Y'$, and also on the $\epsilon$-ball $B'$ of
the normal bundle of $\A_3^{\red,\sing}[\ell]$ at each point of $Y-Y'$. Set
$\Dtilde = B\cup B'$ and $\Wtilde = \partial \Dtilde$. Denote their
pushforwards to $\A_3$ by $D$ and $W$. Then $W$ is a 5-cycle in $\M_3$ which is
generically an $S^3$-bundle over $X'$. Note that $D$ is a 6-chain in $\A_3$
with $\partial D = W$.

\begin{figure}[!ht]
\epsfig{file=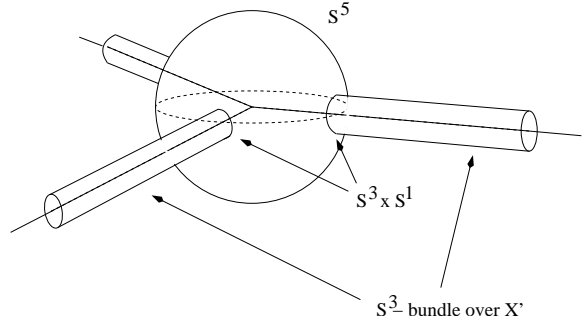, width=3.0in}
\caption{Picture of $W$ near a triple point}
\end{figure}

On the other hand, Looijenga's computation of the rational homology of $\M_3$,
implies that $W$ bounds a rational 6-chain $E$ in $\M_3$. Set $Z = D - E$.
This is a rational 6-cycle in $\A_3$. By construction, we have:

\begin{proposition}
The cycle $Z$ intersections $\A_3^{\red,\sing}$ transversally (in the orbifold
sense), and the intersection number of $Z$ with $\A_3^{\red,\sing}$ is non-zero.
\qed
\end{proposition}

\begin{corollary}
The class of $Z$ is non-trivial in $H_6(\A_3;\Q)$ and the class of
$\A_3^{\red,\sing}$ is non-trivial in $H^6(\A_3;\Q)$. \qed
\end{corollary}

The proof of the non-triviality of $\lambda^3$ is completed by the following
result.

\begin{proposition}
The class of $\A_3^{\red,\sing}$ in $H^6(\A_3;\Q)$ is a non-zero multiple of
$\lambda^3$.
\end{proposition}

\begin{proof}

Denote the closure in $\Abar_3$ of a subvariety $X$ of $\A_3$ by $\Xbar$.
Set $\partial X = \Xbar - X$. Since $\A_3$ is a rational homology manifold,
the sequence
$$
H_8(\Abar_3^\red,\partial \A_3^\red;\Q) \to H^4(\A_3;\Q) \to H^4(\M_3;\Q)
$$
is exact. Since $\A_3^\red$ is irreducible of dimension 4, the left hand
group is one-dimensional and spanned by the fundamental class of $\A_3^\red$.
Since the middle group is one-dimensional and the right hand group trivial,
we see that the class of $\A_3^\red$ spans $H^4(\A_3;\Q)$. On the other hand,
since $\lambda$ is ample, and $\A_3$ contains a complete surface, $\lambda^2$
also spans $H^4(\A_3;\Q)$. It follows that there is a non-zero rational number
$c$ such that
$$
\lambda^2 = c [\A_3^\red] \text{ in } H^4(\A_3;\Q).
$$

Denote the determinant of the Hodge bundle by $\L$. The class $\lambda^3$
is represented by the divisor of a section of the restriction of $\L$
to $\A_3^\red$. This can be computed by pulling back along the mapping
$\A_1 \times \A_2 \to \A_3^\red$.
Since the Picard group of $\A_1$ is torsion, we see that the pullback of $\L$
is represented mod-torsion by $\A_1 \times D$, where $D$ is a cycle
representing $\lambda$ on $\A_2$. But the cusp form $\chi_{10}$ of weight 10 on
$\A_2$ is a section of $\L^{10}$, and has divisor supported on $\A_2^\red$, we
see that $\lambda$ is represented by a non-zero rational multiple of the cycle
$\A_1\times \A_2^\red$ on $\A_1 \times \A_2$, and by a non-zero multiple of
$\A_3^{\red,\sing}$ in $H^4(\A_3^\red;\Q)$. This implies that $\lambda^3$
is represented by a non-zero rational multiple of $\A_3^{\red,\sing}$ in
$H^4(\A_3;\Q)$.

\end{proof}

\section{The rational cohomology of $\Abar_2$ and $\Abar_3$}

Note that $\Abar_0$ is just a point. Suppose now that $g > 0$. Then
$$
\Abar_g = \A_g \amalg \Abar_{g-1}.
$$
Choose a triangulation of $\Abar_g$ such that $\Abar_{g-1}$ is a subcomplex.
Let $U_g$ be a regular PL neighbourhood of $\Abar_{g-1}$ in $\Abar_g$. Set
$U_g^\ast = U_g - \Abar_{g-1}$. This has the homotopy type of $\partial U_g$,
which is a rational homology manifold (as $\Abar_g[3]$ is smooth) of dimension
$2d_g -1$, where $d_g = g(g+1)/2$ is the dimension of $\A_g$. The cohomology
of $U_g^\ast$ has a mixed Hodge structure and the cup product
$$
H^{k-1}(U_g^\ast;\Q) \otimes H^{2d_g - k}(U_g^\ast;\Q)
\to H^{2d_g -1}(U_g^\ast;\Q) \cong \Q(-d_g)
$$
is a perfect pairing of mixed Hodge structures (MHSs) (see \cite{durfee-hain},
for example).

Since $\A_g$ is a rational homology manifold, Lefschetz duality gives an
isomorphism
$$
H_c^k(\A_g,\Q) \cong \Hom(H^{2d_g-k}(\A_g),\Q(-d_g)).
$$
The standard long exact sequence
\begin{equation}
\label{gysin}
\cdots \to H^{k-1}(U_g^\ast) \to H_c^k(\A_g;\Q) \to H^k(\A_g;\Q)
\to H^k(U_g^\ast;\Q) \to \cdots
\end{equation}
is exact in the category of MHSs. These facts will allow us to compute the
cohomology of $U_2^\ast$ and $U_3^\ast$.

The second step in the computation will be to use the Mayer-Vietoris sequence
associated to the covering $\Abar_g = \A_g \cup U_g$
\begin{equation}
\label{mayer-viet}
\cdots \to H^k(\Abar_g) \to H^k(\A_g) \oplus H^k(\Abar_{g-1}) \to H^k(U_g^\ast)
\to H^{k+1}(\Abar_g) \to \cdots
\end{equation}
associated to the covering $\Abar_g = \A_g \cup U_g$, which is exact in the
category of MHS, to compute the cohomology of $\Abar_2$ and $\Abar_3$.

For determining the ring structure and also for seeing that some maps in
these long exact sequences are non-trivial, it is useful to note that since
$\lambda$ is the class of an ample line bundle on $\Abar_g$, the rational
cohomology ring of $\Abar_g$ contains the ring $\Q[\lambda]/(\lambda^{d_g+1})$.

\subsection{Proof of Proposition~\ref{satake_2}}
Since $H^\dot(\A_2;\Q)$ is $\Q$ in degree 0, $\Q(-1)$ in degree 2, and 0
otherwise, $H_c^\dot(\A_2;\Q)$ is $\Q(-2)$ in degree 4, $\Q(-3)$ in degree
6, and 0 otherwise. Using the sequence (\ref{gysin}) and the fact
that the cohomology of $U_2^\ast$ satisfies Poincar\'e duality, we have
$$
H^j(U_2^\ast;\Q) \cong
\begin{cases}
\Q & j = 0; \cr
\Q(-1) & j = 2; \cr
\Q(-2) & j = 3; \cr
\Q(-3) & j = 5.
\end{cases}
$$
Putting this into the Mayer-Vietoris sequence (\ref{mayer-viet}), and using
the fact that $\Abar_1$ is $\P^1$, we obtain the result.

\subsection{Proof of Proposition~\ref{satake_3}}

By duality, $H_c^\dot(\A_3;\Q)$ is an extension of $\Q(-3)$ by $\Q(0)$ in
degree 6, $\Q(-4)$ in degree 8, $\Q(-5)$ in degree 10, and $\Q(-6)$ in degree
12. Using the sequence (\ref{gysin}) and the facts that $H^j(\A_3;\Q)$ vanishes
when $j \ge 7$ and $H_c^j(\A_3;\Q)$ vanishes when $j\le 5$, we have
$$
H^j(U_3^\ast;\Q) \cong
\begin{cases}
H^j(\A_3;\Q) & j < 5;\cr
H_c^{j+1}(\A_3;\Q) & j \ge 7.
\end{cases}
$$
In degrees 5 and 6 we have the exact sequence
$$
0 \to H^5(U_3^\ast;\Q) \to H_c^6(\A_3;\Q) \stackrel{\alpha}{\to} H^6(\A_3;\Q)
\to H^6(U_3^\ast;\Q) \to 0.
$$

It follows from \cite[8.2.2]{deligne} that the image of $\alpha$ is all of
$W_6 H^6(\A_3;\Q)$, so that $\alpha$ is non-zero. Since the sequence is exact
in the category of MHSs, it follows that $H^5(U_3^\ast;\Q)$ is $\Q(0)$ and
$H^6(U_3^\ast;\Q)$ is $\Q(-6)$.

Since $\Abar_3$ is projective and $\lambda$ is the class of a projective
imbedding, the restriction of $\lambda^j$ to $U_3$ is non-zero when $j=1,2,3$.
It follows from this and the computations above that $\lambda$ and $\lambda^2$
restrict to non-trivial classes in the rational cohomology of $U_3^\ast$.

Putting all of this into the Mayer-Vietoris sequence (\ref{mayer-viet}),
easily gives the computation of $H^j(\Abar_3;\Q)$ when $j\neq 6$. The
computation of $H^6(\Abar_3;\Q)$ follows as the sequence
$$
0 \to H^5(U_3^\ast;\Q) \to H^6(\Abar_3;\Q)
\to H^6(\A_3;\Q)\oplus H^6(\Abar_2;\Q)
\to H^6(U_3^\ast;\Q) \to 0
$$
is exact in the category of MHSs.

\end{document}